\documentclass[10pt]{article}
\usepackage[tbtags]{amsmath}
\usepackage{epsfig,amstext,amssymb,amsthm,latexsym}
\usepackage{amssymb,color}

\definecolor{c20}{rgb}{0.,0.7,0.}
\definecolor{c30}{rgb}{0.,0.,1.}
\definecolor{c40}{rgb}{1,0.1,0.7}
\definecolor{c50}{rgb}{1,0,0}

\catcode`\@=11 \@addtoreset{equation}{section} \catcode`\@=12

\allowdisplaybreaks[4] 
\newcommand{\nwc}{\newcommand}
\nwc{\COM}[1]{}

\nwc{\vs}[1]{\vskip #1 cm}

\newtheorem{theo}{Theorem}[section]
\newtheorem{sat}[theo]{Proposition}
\newtheorem{de}[theo]{Definition}
\newtheorem{lem}[theo]{Lemma}

\newtheorem{korr}[theo]{Corollary}

\newtheorem{exxa}[theo]{Example}

\newcommand{\nelem}[1]{{Lemma \ref{#1}}}
\newcommand{\neprop}[1]{{Proposition \ref{#1}}}
\newcommand{\netheo}[1]{{Proposition \ref{#1}}}

\newcommand{\kb}[1]{\boldsymbol{#1}}
\newcommand{\vk}[1]{\kb{#1}}

\def\FRE{\mbox{Fr\'{e}chet }}

\def\X{\vk{X}}

\newcommand{\ve}{\varepsilon}

\newcommand{\abs}[1]{\lvert #1 \rvert}
\newcommand{\Abs}[1]{ \Bigl \lvert #1 \Bigr \rvert}

\newcommand{\E}[1]{\mathbb{E}\{#1\}}

\newcommand{\pk}[1]{\mathbb{P} \{#1\} }

\newcommand{\R}{\mathbb{R}}

\newcommand{\N}{\mathbb{N}}
\newcommand{\inr}{\in \R}

\newcommand{\inn}{\in \N}

\newcommand{\ldot}{,\ldots,}

\newcommand{\limit}[1]{\lim_{#1 \to   \infty}}

\newcommand{\todis}{\stackrel{d}{\to}}

\newcommand{\equaldis}{\stackrel{d}{=}}

\newcommand{\BQN}{\begin{eqnarray}}
\newcommand{\EQN}{\end{eqnarray}}
\newcommand{\BQNY}{\begin{eqnarray*}}
\newcommand{\EQNY}{\end{eqnarray*}}
\newcommand{\BS}{\begin{sat}}
\newcommand{\ES}{\end{sat}}
\newcommand{\BL}{\begin{lem}}
\newcommand{\EL}{\end{lem}}
\newcommand{\BT}{\begin{sat}}
\newcommand{\ET}{\end{sat}}
\newcommand{\BK}{\begin{korr}}
\newcommand{\EK}{\end{korr}}
\newcommand{\BD}{\begin{de}}
\newcommand{\ED}{\end{de}}
\newcommand{\BIT}{\begin{itemize}}
\newcommand{\EIT}{\end{itemize}}
\newcommand{\BDI}{\begin{description}}
\newcommand{\EDI}{\end{description}}

\newcommand{\BEX}{\begin{exxa}}
\newcommand{\EEX}{\end{exxa}}

\newcommand{\QED}{\hfill $\Box$}
\newcommand{\IF}{\infty}
\def\kal#1{{\cal{ #1}}}
\def\fracl#1#2{\biggr( \frac{#1}{#2} \biggl) }

\newcommand{\prooftheo}[1]{ \textsc{Proof of Proposition} \ref{#1} }

\newcommand{\prooflem}[1]{\textsc{Proof of Lemma} \ref{#1}}

\def\X{\vk{X}}

\begin{document}

\centerline{\Large Asymptotics of the convex hull of spherical samples}

        \vskip 1.8 cm
        \centerline{\large Enkelejd Hashorva }
        \vskip 0.3 cm

\centerline{\textsl{Department of Actuarial Science, Faculty of Business and Economics}}

        \centerline{\textsl{HEC Lausanne, University of Lausanne}}

 {\bf Abstract:} In this paper we consider the convex hull of a spherically symmetric sample in $\R^d$. Our main contributions are some new asymptotic results for the expectation of the number of vertices, number of facets, area and the volume of the convex hull assuming that the marginal distributions are in the Gumbel max-domain of attraction. Further, we briefly discuss two other models assuming that the marginal distributions are regularly varying or $O$-regularly varying.



\section{Introduction}
Let $\vk{X}_1 \ldot \vk{X}_n, n\ge 2$ be independent random vectors in $\R^d, d\ge 2$ and denote by $CH[\X_1 \ldot \X_n]$ their convex hull. Distributional and asymptotical properties of the random polytope $CH[\X_1 \ldot \X_n]$ are discussed by many authors, see e.g., R\'{e}ny and Sulanake (1963), Efron (1965),  Raynaud (1970), Carnal (1970),
Eddy and Gale (1981), Groeneboom (1988),  Aldous et al.\ (1991), Carnal and H\"usler (1991), Dwyer (1991), Hueter (1992, 1999, 2004, 2005),
Reitzner (2002, 2004), Buchta (2005), B\'{a}r\'{a}ny and Vu (2007),  Mayer and Molchanov (2007) 
 and the references therein.

In this paper we deal with spherically  symmetric  random vectors assuming the stochastic representation
\BQN\label{eq:sph}
\X_i \equaldis R \vk{U}, \quad i=1 \ldot n,
\EQN
with $R>0$ almost surely being independent of $\vk{U}$ which is uniformly distributed on the unit hypersphere of $\R^d$  (here $\equaldis$ and below $^\top$ stand for the equality of the distribution function, and the transpose sign, respectively).

Next, if $\X=(X_1 \ldot X_d)^\top$ is a spherically symmetric random vector with stochastic representation \eqref{eq:sph}, then in view of Cambanis et al.\ (1981)
\BQN\label{eq:camb}
X_i^2 &\equaldis & R^2 B_{1/2, (d-1)/2}, \quad i=1 \ldot d,
 \EQN
with $B_{1/2,(d-1)/2}$ a Beta distributed random variable with parameters $1/2,(d-1)/2$ being independent of $R$.\\
Since $X_k, k\le d$ are symmetric about $0$ by \eqref{eq:camb} $X_k,k\le d$ have the same distribution function denoted by $Q_d$.

Our main interest lies in the asymptotic properties of $CH[\X_1 \ldot \X_n]$; specifically we focus on the asymptotic behaviour of the expectation of the number of the vertices, facets, the surface area and the 
volume of the convex hull. Interesting asymptotic results for these quantities are derived in the seminal paper Carnal (1970) under explicit assumptions on the tail asymptotics of the distribution function $F$ of $R$ (bivariate setup $d=2$).\\
\def\kG{\kal{G}}
In fact, from the extreme value point of view, Carnal assumed that  $F$ is in the max-domain of attraction (MDA) of a univariate extreme value distribution $\kG$. It is well-known that $\kG$ is either the Gumbel distribution
$\Lambda(x)= \exp(- \exp(-x)), x\inr$, the \FRE distribution $\Phi_\gamma(x)=\exp(-x^{-\gamma}),x>0,\gamma>0 $, or
the Weibull distribution $\Psi_\gamma(x)= \exp(-\abs{x}^\gamma), x< 0$.
Naturally, we raise the question whether Carnal's results can be derived under asymptotic restrictions on $Q_d$?
The answer is positive when $Q_d$ is in the MDA of some univariate distribution function, see Section 3.

Dwyer (1991) extends Carnal's finding to the multidimensional setup assuming again that $F$ is in the MDA of $\kG$.
In the latter paper it is demonstrated that the investigation of the expectation of the number of vertices and facets is of interest for
determine the running time of algorithms for constructing a representation of the facial lattice of the convex
hull of a given point set.

In the Gumbel case ($\kG=\Lambda$) the results of Carnal (1970) and Dwyer (1991) are valid for special distribution functions $F$
with light exponential tails and infinite upper endpoint (referred below as Carnal distributions).
Asymptotic results for the expectation of the number of the vertices of the convex hull
are to date not available when $F$ is in the Gumbel MDA  and has a finite upper endpoint.

Without going to mathematical details, we briefly mention the main contributions of this paper:

a) Making use of extreme value theory, we extend the known results for the Carnal distributions  $F$ to the larger class of univariate distribution functions in the Gumbel MDA.
Furthermore, we obtain asymptotic estimates for the expectation and the variance of the number of vertices of the convex hull as well as a CLT extending a fundamental theorem of Hueter (1999).

b) We show that several existing results can be derived with similar assumptions on the marginal distribution function $Q_d$ giving a positive answer to the above question.

c) A new result derived in this paper is the boundedness  of the sequence of the expectation of the number of the vertices of the convex hull if either $\overline{F},$ or  $\overline{Q}_d$  are $O$-regularly varying.

Organisation of the paper: The main results are presented in Section 3 followed by a section dedicated to the proofs. We conclude the paper with an Appendix.

\section{Preliminaries}
We introduce first our notation, provide few results from extreme value theory, and review some known results for the convex hull $CH[\X_1 \ldot \X_n]$ of a spherically  symmetric sample $\X_1 \ldot \X_n$ as given in the Introduction.

If $H$ is the distribution function of a random variable $Y$ (henceforth  abbreviated as $Y\sim H$), then we write $\overline{H}:= 1- H$ for its survival function. Further we define
the generalised inverse of $H$ by $H^{-1}(s):=\inf\{x: H(x)\ge s\}$ and denote by $x_H:=\sup\{x:H(x)< 1\}$ the upper endpoint of $H$. We use similar notation for other distributions.\\
Throughout in the following $B_{\alpha, \beta}$ stands for a Beta random variable with positive parameters $\alpha, \beta$ with density function
$$x^{\alpha-1}(1- x)^{\beta-1}\frac{\Gamma(\alpha+\beta)}{\Gamma(\alpha)\Gamma(\beta)}, x\in (0,1),$$
where $\Gamma(\cdot)$ is the Euler Gamma function.

From extreme value theory (see e.g., Reiss (1989), Embrechts et al.\ (1997), Falk et al.\ (2004)) the univariate distribution function $N$ is in the MDA of the univariate distribution function $\kG$, if for some constants $a_n>0, b_n, n\inn$
\BQN\label{eq:maxH}
 \limit{n} \sup_{t\inr} \Abs{N^n(a_nt+ b_n)- \kG(t)} &=&0.
\EQN
As mentioned above only three choices are possible for $\kG$, namely $\Lambda, \Phi_\gamma$ or $\Psi_\gamma$, with$ \gamma \in (0,\IF)$. When $\kG=\Lambda$ the upper endpoint $x_N$ of $N$ can be finite of infinite. For both other cases, $x_N$ is either finite (Weibull) or infinite (Fr\'{e}chet). The characterisation of both Weibull and \FRE max-domain of attractions is closely related to the concept of the regularly varying functions. In the following a positive measurable function $\kal{L}$ is called slowly varying at infinity if $\limit{u} \kal{L}(u s)/\kal{L}(u)=1$ for any $s>0$. A regularly varying function with index $\gamma\inr$ is the product of some $\kal{L}(x)$ with $x^\gamma$.

Next, we briefly review some known results for the convex hull. Let $v_n, f_n$ denote the number of the vertices and the facets of $CH[\X_1 \ldot \X_n]$. Referring to Dwyer (1991) we may write
\BQN\label{eq:vn1}
 \int_0^\IF Q^{n-1}_d(s) \abs{d \overline{F}(s)}
&\le & \frac{\E{v_n}}{n}
\le  (1+ o(1))2^{d-1} \int_0^\IF [1 - 2^{-d+1}\overline{Q}_d(s)]^{n-1} \abs{d \overline{F}(s)}, \quad n\to \IF
\EQN
(we abuse slightly the notation writing $P^k(s)$ instead of $(P(s))^k, k\inr$ for $P$ some arbitrary function).  Furthermore
\BQN\label{eq:vn2}
\E{v_n} \le \E{f_n} \sim \frac{n^d}{d} \kappa_d \int_0^\IF \delta_1(r) q^d_d(r) \exp(- n \overline{Q}_d(r))\,d r, \quad n\to \IF,
\EQN
where $q_d$ is the density function of $Q_d$ (which exists, see \eqref{PKS:14:2} in Appendix) and $\delta_1(r)$ is bounded by  (see Lemma 4 in Dwyer (1991))
\BQN
\delta_1(r) &\le&  (1+o(1)) \frac{d \tau_{d-1} \kappa_{d-1} }{q_d(r)} \int_r^\IF(u^2- r^2)^{d-2} u \, dF(u), \quad r\to \IF,
\EQN
where
\BQN \label{eq:mu:tau}
\tau_d:= \sqrt{d}\frac{(1+1/d)^{(d+1)/2}}{\Gamma(d+1)} , \quad \kappa_d:=\frac{2 \pi ^{d/2}}{\Gamma(d/2)}.
\EQN
In \eqref{eq:vn2} and below $a_u \sim b_u, u \uparrow \omega ,$ with $\omega\in (-\IF, \IF]$ means that $\lim_{u \uparrow \omega} a_u/b_u=1$. Further, we write $a_n \sim b_n$ instead of $\limit{n} a_n/b_n=1$.\\
 In the two-dimensional setup $d=2$ (see Carnal (1970))
\BQN\label{eq:vn1:2}
2\E{v_n}&\sim & n^2 \int_0^\IF Q^{n-2}_2(s)\abs{d \overline{H}(s)},
\EQN
with 
\BQN \label{A1}
H&\thicksim& \min(R_1,R_2) \sqrt{B_{1/2,1/2}}, \quad R_1\equaldis R_2\equaldis R,
 \EQN
where $R_1,R_2, B_{1/2,1/2}$ are mutually independent, see {\bf A1} for the proof.

Referring again to Dwyer (1991), we have for the surface area $A_n$ and the volume $V_n$ of the convex hull
\BQN\label{eq:An:d}
\E{A_n}& \sim & \frac{n^d}{d} \tau_d \int_0^\IF \delta_2(r) q^d_d(r)\exp(- n \overline{Q}_d(r))\, dr
\EQN
and
\BQN\label{eq:Vn:d}
\E{V_n}& \sim & n^d \tau_d \int_0^\IF r \delta_2(r) q^d_d(r)\exp(- n \overline{Q}_d(r))\, dr,
\EQN
where  $\delta_2$ can be bounded asymptotically by
\BQN\label{eq:Delta2}
\delta_2(r) &\le &
(1+o(1))\frac{ d \tau_{d-1} \kappa_{d-1}}{q_d(r)} \int_r^\IF (u^2- r^2)^{(3d-5)/2} u \, dF(u), \quad r\to \IF.
\EQN
In the bivariate setup
\BQN\label{eq:An:d}
2\E{A_n}&\sim &   n^2\int_0^\IF Q^{n-2}_d(s) \abs{d \overline{K}(s)},
\EQN
with
$$\overline{K}(s)= 1- K(s)=\frac{1}{\pi} \biggl [\int_s^\IF \sqrt{y^2- s^2} \, d F(y) \biggr]^2, \quad s\ge 0.$$
We note in passing that $K(s), s\ge 0$ is  continuous, see Appendix {\bf A2}.

\section{Main Results}
Asymptotic results for the sequence of the expectation of the number of vertices $\E{v_n},n\ge 1$ can be obtained (when $d=2$) by investigating the asymptotic behaviour $(n\to \IF$) of $\int_{\R} Q_2^n(s)\, \abs{d \overline{H}(s)}$. Our \nelem{lem:2} turns out to be quite useful;  furthermore it sheds some light explaining the role of extreme value theory in our analysis. More specifically, in view of \nelem{lem:2} $\E{v_n},n\ge 1$ is regularly varying sequence $(n \to \IF)$ if and only if (iff)  the tails of $Q_2$ and $F$ satisfy a certain asymptotic condition (see \eqref{lem2:2}). In particular
\BQN\label{eq:c}
\limit{n} \E{v_n}&=& c\in (0,\IF)
\EQN
iff
\BQN\label{eq:de}
\lim_{u \uparrow x_H}\frac{\overline{H}(u)}{\overline{Q}_2^2(u)}&= & c \in (0,\IF).
\EQN
On the other hand, $H$ and $Q_2$ are in a strong relation with the distribution function $F$ via \eqref{eq:camb}. Thus
it is not straightforward to check whether \eqref{eq:de} holds for some given $F$. One simple
instance is when  for some positive constants $c_1,c_2\in (0,\IF)$
\BQN \label{eq:const}
\overline{H}(u)&\sim &c_1 \overline{F}^2(u), \quad \overline{Q}_2(u)\sim c_2 \overline{F}(u),\quad u \uparrow x_H
\EQN
implying that \eqref{eq:c} is valid with $c:=c_1/c_2^2$. \\

If $F$ is in the MDA of an extreme value distribution, then by Berman (1992) it follows that both $Q_2$ and $H$ are in the same MDA, and further the asymptotics of $\overline{Q}_2(u)$ and $ \overline{H}(u)$ as $u\to x_F$ are
determined by $\overline{F}(u)$ and some known functions.\\
In view of Hashorva and Pakes (2010) also the converse holds, i.e.,  $F$ is in the MDA of an extreme value distribution function iff  $Q_d,d\ge 2$ is in the MDA of the same extreme value distribution function. Consequently, all known results with  $F$ in the MDA  of some univariate extreme value distribution can be retrieved if we impose instead the same assumption on $Q_d$.\\
We deal first with the Gumbel case; when $x_F\in (0,\IF)$ no asymptotic results for the quantities of interest are known to date. When $x_F=\IF$ we have both results of Carnal and Dwyer for any $F$ being a Carnal distribution function. We conclude this section by briefly discussing both the \FRE and Weibull max-domains of attraction.\\

\subsection{Gumbel Tails}
It is well-known that condition \eqref{eq:maxH} is valid for some univariate distribution function $N$ with $\kG=\Lambda$ and upper endpoint $x_N\in (-\IF, \IF]$, iff for some positive scaling function $w$
\BQN
\label{eq:rdfd}
\lim_{u \uparrow x_N} \frac{\overline{N}(u+x/w(u))} {\overline{N}(u)} &=& \exp(-x),\quad \forall x\inr.
\EQN
Furthermore
\BQN \label{eq:wN}
 w(u)&\sim &\frac{\overline{N}(u)}{\int_u^\IF \overline{N}(s)\, ds}, \quad u \uparrow x_N.
 \EQN
So far in the literature the Gumbel MDA  assumption on $F$ has not been explicitly assumed.
An elegant simplification of this assumption is suggested in Carnal (1970) which has been used in several following papers (Eddy and Gale (1981), Dwyer (1991), Hueter (1999, 2005, 2005)).
More specifically, Carnal (1970) considers distribution functions $F$ satisfying (for all large $x$)
\BQN\label{eq:carn}
 x&=&\kal{L}(1/\overline{F}(x)),
 \EQN
where $\kal{L}$ is a monotone increasing slowly varying function at infinity.
We refer to \eqref{eq:carn} as the Carnal tail condition and to such  $F$ as Carnal distributions. As shown in Carnal (1970) if $\kal{L}(s)=\exp(\int_1^s \ve(s)/s\, ds),s\ge 0$ with  $\limit{s}\ve(s)=0$, then  Carnal distributions have the representation
$$\overline{F}(x)=\exp(- \int_0 ^x 1/(\eta(s)s)\, ds), \quad x\ge 0,$$
with $\eta(s)=\ve(1/\overline{F}(s)),s>0$. In the aforementioned paper (see also Dwyer (1991), Hueter (1999)) the function $\eta$ satisfies some smoothness conditions being further positive and monotone non-decreasing.\\
If $N\in GMDA(w)$ we define next
\BQNY
\xi_N(n)&:=&b_n/a_n, \quad b_n:=N^{-1}(1- 1/n), \quad a_n:=1/ w(b_n), \quad n> 1.
\EQNY
The constants $a_n,b_n$ are such that \eqref{eq:maxH} holds with $\kG=\Lambda$. Further, it is well-known  (see e.g., Resnick (2008)) that both $N^{-1}(1- 1/n)$ and $\xi_N(n)$ are slowly varying functions at infinity. As will be shown next this fact,  \nelem{lem:2} and \netheo{theo:0} (see Appendix) are the key ingredients needed to derive the tail asymptotics of the quantities of interest.\\


\BT \label{theo:m} Let $F,H,K,Q_d,d\ge 2$ be as in the previous section, and let $CH[\X_1 \ldot \X_n]$ be the convex hull of the random points $\X_1\ldot \X_n$ which are  independent with stochastic representation \eqref{eq:sph}, where $R\sim F$. If $Q_d \in GMDA(w)$ or $F\in GMDA(w)$ we have:\\
a) As $n\to \IF$
\BQN
 Q^{-1}_d(1-1/n)\sim F^{-1}(1-1/n), \quad  \text{ and } \quad  \xi_{Q_d}(n) \sim \xi_F(n).
 \EQN
b) If $d=2$, then
\BQN\label{eq:vn1:res:1}
\E{v_n}& \sim & 2 \sqrt{\pi \xi_{Q_2}(n)},
\EQN
and for $d \ge 2$ and some $\ve\in (0,\IF)$
\BQN\label{eq:vn1:res:1}
(1- \ve) \xi_{Q_d}(n)^{(d-1)/2}\le 2^{(d-3)/2}\frac{\Gamma(d/2)} {\sqrt{\pi}}\E{v_n}& \le (1+ \ve)  4^{d-1}\xi_{Q_d}(n)^{(d-1)/2}.
\EQN
c) If $d=2$, then
\BQN \label{eq:vn1:res:2}
\E{A_n}&\sim & \pi [Q^{-1}_2(1-1/n)]^2, 
\EQN
and for $d\ge 2$
\BQN\label{eq:vn1:res:3}
\E{A_n}&\le & (1+o(1)) \Gamma(d+1)\fracl{4d \sqrt{\pi}}{d-1}^{d-1} (Q_d^{-1}(1-1/n))^{d},  \quad n\to \IF.
\EQN
d) For any $d\ge 2$
\BQN\label{fn}
\E{f_n} &\le &(1+o(1)) \sqrt{d} \fracl{ 8 \pi d}{d-1}^{(d-1)/2} (\xi_{Q_d}(n))^{-(d-1)/2}, \quad n\to \IF
\EQN
and
\BQN\label{vn}
\E{V_n}  &\le &(1+o(1)) \Gamma(d) \fracl{ 4 d\sqrt{\pi} }{d-1}^{d-1} (Q^{-1}_d(1-1/n))^{d}, \quad n\to \IF
\EQN
are valid.
\ET

{\bf Remarks:} {\it (a) In Carnal's notation $L(n)=F^{-1}(1-1/n)$ and $\ve(n):=1/\xi_{F}(n), n\ge 1$. If the upper endpoint of $F$ is finite, say $x_F=1$, then clearly $\limit{n}L(n)=1$. \\
(b) A misprint appears in the upper bound for $\E{v_n}$ in Dwyer (1991), p.126. The upper bound therein should be multiplied by $2^{d-1}$, see \eqref{eq:vn1:res:1} above.\\
(c) For any $N\in GMDA(w)$ with scaling function $w$ defined by \eqref{eq:wN} we have  (see e.g., Resnick (2008))
\BQN\label{tn2}
u w(u) \to \IF, \quad \text{ and if $x_N< \IF$   } w(u)(x_N- u)&=& \IF, \quad u\uparrow x_N.
\EQN
Consequently, \eqref{eq:vn1:res:1} implies $\limit{n}\E{v_n}=\IF$.\\
(c) Utilising the expression (1.7) which gives an asymptotic formula for $\E{l_n}$ with $l_n$ the perimeter of the convex hull (d=2), it follows
that when $F$ or $Q_d$ are in the Gumbel MDA with some scaling function $w$, then we have
\BQN
\E{l_n} &\sim &2 \pi Q_2^{-1}(1-1/n).
\EQN
}

\neprop{theo:m} provides asymptotic upper and lower bounds for $\E{v_n}$.\\
 Hueter (1999) was able to give the exact asymptotic behaviour of the first and the second moment of $v_n$; moreover a key central limit theorem  was derived therein
 by developing Groeneboom's technique (see Groeneboom (1988)) in higher dimensions.

Next we extend Hueter's CLT theorem which has been shown for Carnal distributions
by considering  a general spherical random vector with marginal distribution or distribution
of the associated random radius in the Gumbel max-domain of attraction.\\

\BT \label{theo:Hu}
Let $v_n$ denote the number of the vertices of $CH[\X_1 \ldot \X_n]$, where $\X_i, i\ge n$ are independent with stochastic representation \eqref{eq:sph}. Suppose that $F(0)=0, x_F=\IF$ and set
$a_n:=1/w(b_n), b_n:=Q_d^{-1}(1-1/n),n>1$ where $Q_d$ and $F$ are related by \eqref{eq:camb}.
 If either $F\in GMDA(w)$ or $Q_d\in GMDA(w)$, and $a_n b_n \not \to \IF$ as $n\to \IF$, then we have the convergence in distribution
\BQN
\frac{v_n- \lambda_d ( b_n/a_n)^{(d-1)/2}}{\sqrt{ {\bf Var} \{ v_n \} }} &\todis &Z , \quad n\to \IF,
\EQN
with ${\bf Var} \{ v_n \} \sim \lambda^*_d (b_n/a_n)^{(d-1)/2}, \lambda_d,\lambda_d^*\in (0,\IF)$, and $Z$ a standard Gaussian random variable.
\ET

{\bf Remarks}:  {\it (a) The above proposition gives also the asymptotics of $\E{v_n}$ and ${\bf Var}\{v_n\}, n\to \IF$. It turns out that the expectation and the variance of the number of the vertices of the convex hull differ by a constant, and are both slowly varying functions.\\
b) The condition $a_n b_n \not \to \IF$ as $n\to \IF$ implies a certain asymptotic behaviour of the density function $q_2$ of $Q_2$. More precisely, in view of \netheo{theo:m}
\BQNY
q_2(u)&\sim& w(u) \overline{Q}_2(u)=(1+o(1)) \frac{1}{2 \pi} \fracl{ w(u)}{u}^{1/2} \overline{Q}_2(u),
\quad u\to \IF,
\EQNY
hence
\BQNY
q_2(b_n)&\sim & \frac{1}{2 \pi} \fracl{ 1}{a_n b_n}^{1/2} \overline{Q}_2(b_n)\sim \frac{n}{\sqrt{a_n b_n}},
\EQNY
implying  $q_2(b_n)/n \not \to 0$ as $n\to \IF$.\\
(c) In \netheo{theo:Hu} if $F$ has a finite upper endpoint $x_{F}\in (0,\IF)$, then $\limit{n} b_n=x_{F}$ and $\limit{n} a_n=0$, hence $\limit{n}a_n b_n=0$.\\
We conjecture that \netheo{theo:Hu}, and in particular the asymptotics of $\E{v_n}$ and ${\bf Var}\{v_n\}$, are
also valid with the same constants (not depending on $F$), when $F\in GMDA(w)$ with $x_{F} \in (0,\IF)$. In the 2-dimensional setup this is true for the asymptotics of $\E{v_n},n\to \IF$ (see \eqref{eq:vn1:res:1} above).
} \\

We give next two illustrating examples.

{\bf Example 1.} Let $\X_i \equaldis R  \vk{U}, i=1 \ldot n $ be independent spherically  symmetric random vectors in $\R^2$. Assume that the distribution function $F$ of the positive random variable $R$ has upper endpoint 1
 satisfying
$$ \overline{F}(u) \sim  a \exp(- b/(1- u)) , \quad u \uparrow 1,$$
with $a,b\in (0,\IF)$. Set $w(u):= b/(1- u)^2, u\in (0,1)$. Since for any $s\inr$
$$ \frac{ \overline{F}(u+ s/w(u))}{ \overline{F}(u)}= (1+o(1))\exp( - b [ 1/(1- u+ s/w(u)) - 1/(1- u)]) \to \exp(-s), \quad u \uparrow 1,$$
then $F\in GMDA(w)$. Further, we have
$$ F^{-1}(1-1/n)\sim 1- b/\ln (a n), \quad w(F^{-1}(1-1/n))=b[\ln (a n)]^2,
\quad n> 1,$$
consequently
$$ \xi_{F}(n)= F^{-1}(1-1/n) w(F^{-1}(1-1/n))\sim b[\ln  n]^2 ,  \quad n\to \IF.$$
Hence in view of \netheo{theo:m} for $d=2$
\BQN\label{new}
 \E{v_n} &\sim&  \sqrt{ 4 b\pi}  \ln n, \quad n\to \IF.
 \EQN
\COM{If $d\ge 3$, then \netheo{theo:Hu} implies
\BQN\label{new:d}
 \E{v_n} &\sim& \lambda_d b^{(d-1)/2} (\ln n)^{d-1}, \quad  {\bf Var}\{v_n\} \sim \lambda_d^*b^{(d-1)/2} (\ln n)^{d-1}, \quad n\to \IF,
 \EQN
with $\lambda_d, \lambda_d^*$ as in the aforementioned proposition. Note in passing that $\lambda_2=2 \sqrt{\pi}$,
and $F$ is not a Carnal distribution.
}
\bigskip
{\bf Example 2.} Under the setup of the previous example, we suppose further that
the marginal distribution function $Q_2$ is in the Gumbel MDA  with scaling function
$$ w(x)= r\theta x^{\theta-1}/(1+ \kal{L}_1(x)), \quad r>0, \theta>0, $$
where $\kal{L}_1$ is a regularly varying function at infinity with index $ \gamma \theta, \gamma<0$. It follows
that
$$\overline{Q}_2(x)\sim \exp(- r x^\theta(1+ \kal{L}_2(x))), \quad x\to \IF,$$
with  $\kal{L}_2$ another regularly varying function at infinity with index $ \gamma \theta$. Consequently, we have
$$b_n:= Q_2^{-1}(1-1/n)\sim  \fracl{\ln n}{r}^{1/\theta}, \quad a_n:= 1/w(b_n)= \frac{b_n^{1- \theta}}{r \theta }, \quad n\to \IF$$
implying
$$ \xi_{Q_2}(n)= \theta \ln n , \quad a_n b_n= \frac{b_n^{2- \theta}}{r \theta }, \quad n>1.$$
Hence, by \netheo{theo:m}
\BQN
\E{v_n}& \sim & 2 \sqrt{\pi \theta \ln n }
\EQN
and $\limit{n} a_n b_n=0$ if $\theta>2$, whereas for $\theta \in (0, 2)$ we have $\limit{n} a_n b_n=\IF$.
Note in passing that if $Q_2$ is the standard Gaussian distribution, then $\theta=1/r=2$ and $\limit{n}a_nb_n= 1$. Further,  we remark that if $\kal{L}_1$ is constant, then $Q_2$ is a Carnal distribution.

\subsection{Regularly and $O$-Regularly Varying Tails}
The survival function $\overline{N}$ is regularly varying (at infinity) with index $\gamma\le  0$ if
\BQN\label{eq:fre1}
\limit{u} \frac{\overline{N}(ux)}{\overline{N}(u)}&=& x^{\gamma}.
\EQN
In view of \netheo{eq:theo:BM1}  (see Appendix {\bf A3}), the survival function $\overline{Q}_d$ satisfies \eqref{eq:fre1} iff $\overline{F}$ also satisfies \eqref{eq:fre1}. Hence the results of Carnal (1970) and Dwyer (1991) can be retrieved assuming the regular variation of $\overline{Q}_d$ instead of that of $\overline{F}$.
As shown in Berman (1992) it is possible to relate the asymptotics of $\overline{F}$ with that of $\overline{Q}_2$, specifically
\BQNY
\overline{Q}_2(u)&\sim & \frac{\Gamma((\gamma+1)/2)}{\sqrt{\pi} \Gamma(\gamma/2+1)}\overline{F}(u), \quad u \to \IF.
\EQNY
Similarly, we find $\overline{F}$ satisfies \eqref{eq:fre1}  iff the survival function $\overline{H}$ is regularly varying with index $2 \gamma$. Moreover as $u\to \IF$
\BQNY
\overline{H}(u)&\sim & \frac{\Gamma(\gamma+1/2)}{\sqrt{\pi} \Gamma(\gamma+1)}\overline{F}^2(u) .
\EQNY
Consequently \eqref{eq:const} implies  that if one of the survival functions $\overline{F}, \overline{Q}_d$ or $\overline{H}$ is regularly varying with index $\gamma \le 0$, then for the bivariate setup ($d=2$) we have
\BQN\label{eq:limfin}
\limit{n} \E{v_n}&=& \frac{\Gamma(\gamma+1/2)[\Gamma(\gamma/2+1)]^2}{[\Gamma((\gamma+1)/2)]^2 \Gamma(\gamma+1)},
\EQN
which is shown in Carnal (1970) assuming that $\overline{F}$ satisfies \eqref{eq:fre1}. Aldous et al.\ (1991) addresses the case that $\overline{F}$ satisfies \eqref{eq:fre1} with $\gamma=0$, which in view of our results is equivalent with $Q_d$ being slowly varying (satisfying \eqref{eq:fre1} with $\gamma=0$).\\
It is interesting that the limit in \eqref{eq:limfin} is finite, thus $\E{v_n},n\ge 1$ is  a bounded sequence.\\
 A natural question that arises is: For what other distribution functions is $\E{v_n},n\ge 1$ a bounded sequence? \\
 We show below that the answer is positive for $\overline{F}$ being a $O$-regularly varying function, meaning that
$$ 0 < \liminf_{u \to \IF} \frac{\overline{F}(ux)}{\overline{F}(u)}\le
\limsup_{u \to \IF} \frac{\overline{F}(ux)}{\overline{F}(u)}< \IF, \quad \forall x> 1.$$

\BT \label{theo:Or} Under the setup of \netheo{theo:m}, if either $\overline{F}$, $\overline{Q}_2,$ or $\overline{H}$ is $O$-regularly varying at infinity, then $\E{v_n},n\ge 1$ is a bounded sequence.
\ET

It is well-known that $F$ is in the \FRE MDA  iff \eqref{eq:fre1} holds for some $\gamma <0$, see e.g., Embrechts et al.\ (1997). When $F$ is in the Weibull MDA
we have a similar behaviour of the survival function at the upper endpoint $x_F$ which is necessarily finite, say $x_F=1$. More specifically
\BQN\label{eq:wei}
\limit{u}\frac{\overline{F}(1- x/u)}{\overline{F}(1- 1/u)}&=&x^\gamma, \quad \gamma >0
\EQN
is equivalent with \eqref{eq:maxH} where $\kG= \Psi_\gamma$.   

Our new results when $F$ is in the Weibull MDA can be derived utilising \neprop{eq:theo:BM1} which implies:\\
Condition \eqref{eq:wei} is equivalent with the fact that the marginal distribution $Q_d$ satisfies \eqref{eq:wei} with $\gamma^*=\gamma+1/2, \gamma\ge 0$. Hence the results of Carnal (1970) for the bivariate setup, and those of Dwyer (1991) for the higher dimensions
hold if either of those conditions are satisfied.\\

\section{Proofs}
\prooftheo{theo:m} We note first that for the proof \netheo{theo:0} presented in the appendix is crucial.\\
a) The claim is easily established if $F$ has a finite upper endpoint, say $x_F=1$. Since also $Q_d$ has the same upper endpoint, the proof follows by the fact that
$$ \limit{n} Q_d^{-1}(1-1/n)=x_{Q_d}=1, \quad \limit{n} F^{-1}(1-1/n)=x_{F}=1.$$
We deal therefore with the case $x_F=\IF$. By Lemma 6.1 in Hashorva (2009) if $N\in GMDA(w)$ and $N$ has an infinite upper endpoint, then
\BQN\label{eq:resn}
 \limit{u} \frac{\overline{N}(c u) (u w(u))^a}{\overline{N}(u)}&=&0
\EQN
for any $c>1$ and $a\inr$. Hence, by \eqref{tn2}, \eqref{eq:resn}, \eqref{eq:loc} and \eqref{theo:0:b}
\BQNY
 Q_d^{-1}(1-1/n) &\sim & F^{-1}(1-1/n), \quad n\to \IF.
 \EQNY
In view of \netheo{theo:0} both $H$ and $Q_d, d\ge 2$ belong to the Gumbel MDA  with the same scaling function $w$, consequently  $\xi_{Q_d}(n) \sim \xi_F(n)$.\\
b) By the assumptions and \netheo{theo:0} we have
\BQNY
\frac{\overline{H}(u)}{ \overline{Q}^2_2(u) \sqrt{\xi_{Q_2}(\overline{Q}_2(u))}} \sim 2 \sqrt{\pi}, \quad u\uparrow x_F.
\EQNY
Furthermore, $H$ and $Q_d,d\ge 2$ are continuous distribution functions. Hence applying \nelem{lem:2} with $\rho=2$ and $l= \Gamma(3) 2 \sqrt{\pi}$ to  \eqref{eq:vn1:2} we obtain
\BQNY
2\E{v_n}&\sim & n^2 \int_0^\IF Q^{n-2}_2(s)\abs{d \overline{H}(s)}
\sim \Gamma(3)2 \sqrt{\pi}.
\EQNY
If $d\ge 2$, then  \eqref{theo:0:a} implies
\BQNY
\frac{ \overline{F}(u)}{\overline{Q}_d(u) \xi_{Q_d}( \overline{Q_d}(u))^{(d-1)/2}} \sim 2^{-(d-3)/2}\frac{\sqrt{\pi}}{\Gamma(d/2)}, \quad u\uparrow  x_F.
\EQNY
Consequently,  by \nelem{lem:2}
\BQNY
n\int_0^\IF Q_d^{n-1}(s) \, d F(s) \sim  2^{-(d-3)/2}\frac{\sqrt{\pi}}{\Gamma(d/2)}\xi_{F}(n)^{(d-1)/2}.
\EQNY
Similarly,
\BQNY
n2^{d-1}\int_0^\IF [1- 2^{-(d-1)} \overline{Q}_d(s)]^{n-1} \, d F(s) \sim
4^{d-1} 2^{-(d-3)/2}\frac{\sqrt{\pi}}{\Gamma(d/2)}  \xi_{F}(n)^{(d-1)/2}
\EQNY
and thus \eqref{eq:vn1:res:1} follows. \\
c) With the same arguments as above we obtain
\BQNY
\frac{ \overline{K}(u)}{  [\overline{Q}(u) Q^{-1}( 1- \overline{Q}(u))]^2} \sim  \pi, \quad u\uparrow x_F,
\EQNY
hence \nelem{lem:2} implies
\BQNY
2\E{A_n}&\sim &   n^2\int_0^\IF Q^{n-2}_2(s) \abs{d \overline{K}(s)} \sim \Gamma(3) \pi [Q^{-1}_2(1-1/n)]^2,
\EQNY
and thus \eqref{eq:vn1:res:2} follows.\\
Next we show the claim for $d\ge 2$.  Applying Lemma 7.6 of  Hashorva (2007) and using \eqref{eq:An:d} 
 we obtain
\BQNY
\int_r^{x_F} (r^2- u^2)^{(3d-5)/2} u \, d F(u)&\sim & \Gamma(3(d -1)/2) \fracl{2r}{w(r)}^{(3d-5)/2} r \overline{F}(r), \quad r\uparrow x_F.
\EQNY
\COM{

\BQNY
\delta_2(r)&\le& (1+o(1)) \frac{d \tau_{d-1} \kappa_{d-1}}{q_d(r)} \Gamma(3(d -1)/2) \fracl{2r}{w(r)}^{(3d-5)/2} r \overline{F}(r)\\
&=&  (1+o(1)) \fracl{2d}{d-1}^{d-1} \frac{ d^2(d-1)!}{2\Gamma((d+1)/2)}  \fracl{r^2}{ w(r)}^{d-1}.
\EQNY
}
Consequently, by \eqref{eq:An:d}, \eqref{theo:0:a} and the fact that $q_d(r) \sim w(r) \overline{Q}_d(r), r\uparrow x_F$ (see \netheo{theo:0})
\BQNY
\E{A_n}&\le& (1+o(1)) \fracl{2d}{d-1}^{d-1} \frac{ d^2(d-1)!}{2\Gamma((d+1)/2)}
\int_0^{x_F} \fracl{r^2}{w(r)}^{d-1} (r w(r))^{d-1}(r) \exp(- \overline{Q}_d(r)) \, dr\\
&\le& (1+o(1)) \fracl{2d}{d-1}^{d-1} \frac{ d^2(d-1)!}{2\Gamma((d+1)/2)}
\int_0^{x_F} r^{2(d-1)} \exp(- \overline{Q}_d(r)) \, dQ_d(r),  \quad r\uparrow x_F.
\EQNY
Since $Q^{-1}(1-1/n)$ is regularly varying as $n\to \IF$ the Abel formula for the Laplace transform yields
\BQNY
\E{A_n}&\le& (1+o(1))\Gamma(d+1)\fracl{4d \sqrt{\pi}}{d-1}^{d-1} (Q_d^{-1}(1-1/n))^{d-1},  \quad n\to \IF,
\EQNY
thus the statement is established. \\
d) By \netheo{theo:0} and Lemma 7.6 in Hashorva (2007) as $r \uparrow x_F$ we obtain
\BQNY
\frac{1}{q_d(r)}\int_r^{x_F} (r^2- u^2)^{d-2} u \, d F(u)&\sim & \Gamma(d-1) \fracl{2r}{w(r)}^{d-2} r \frac{\overline{F}(r)}{q_d(r)}\\
&\sim& \Gamma(d-1) 2^{d-2} r^{d-1} w(r)^{1- d}\frac{w(r)\overline{F}(r)}{q_d(r)}\\
&\sim& \Gamma(d-1) 2^{d-2} 2^{-(d-3)/2} r^{d-1} w(r)^{1- d} (r w(r))^{(d-1)/2} \sqrt{\pi}/\Gamma(d/2) \\
&\sim& \frac{\Gamma(d-1)\sqrt{\pi}}{\Gamma(d/2) } 2^{d/2-1/2} r^{3(d-1)/2} w(r)^{(1- d)/2}.
\EQNY
Hence \eqref{fn} and \netheo{theo:0} implies 
\BQNY
\E{f_n}&\le & (1+o(1)) d \tau_{d-1} \kappa_{d-1}
\frac{\Gamma(d-1)\sqrt{\pi}}{\Gamma(d/2) } 2^{d/2-1/2}  \int_0^{x_F}
r^{3(d-1)/2} w(r)^{(1- d)/2}  q_d(r)^{d-1} \exp(-n\overline{Q}_d(r))\, d Q_d(r)\\
&\le & (1+o(1)) d \tau_{d-1} \kappa_{d-1}
\frac{\Gamma(d-1)\sqrt{\pi}}{\Gamma(d/2) } 2^{d/2-1/2}  \int_0^{x_F}
(r w(r))^{3(d-1)/2}   \overline{Q}_d(r)^{d-1} \exp(-n\overline{Q}_d(r))\, d Q_d(r).
\EQNY
In view of Theorem 4.1 in Hashorva et al.\ (2010) there exists a distribution function $G_d$ such that
\BQNY
\overline{G}_d(u) &\sim  &  d \tau_{d-1} \kappa_{d-1}
\frac{\Gamma(d-1)\sqrt{\pi}}{\Gamma(d/2) } 2^{d/2-1/2}  (u w(u))^{3(d-1)/2}   \overline{Q}_d(u)^{d-1}, \quad u\to \IF.
\EQNY
Applying now \nelem{lem:2} to $ \int_0^\IF \overline{G}_d(r)  \exp(-n\overline{Q}_d(r))\, d Q_d(r)$
establishes \eqref{fn}.

The proof of the last claim follows with similar arguments utilising further \eqref{vn}.  \QED

\prooftheo{theo:Or} Applying \nelem{lem:OR} $\overline{F}$ is $O$-regularly varying if $\overline{Q}_d$ or $\overline{H}$ is $O$-regularly varying, and vice-versa. Next, assume that $d=2$ and $\overline{F}$ is $O$-regularly varying. By \eqref{eq:camb} and \eqref{A1} for any $c>1$ and $u>0$ we have
$$ \pk{B_{1/2,(d-1)/2} > 1/c} \overline{F}(c u) \le\overline{Q}_d(u) \le \overline{F}(u) $$
and
$$ \pk{B_{1/2,(d-1)/2} > 1/c} \overline{F}^2(c u) \le  \overline{H}(u) \le \overline{F}^2(u). $$
Consequently
$$ \fracl{\overline{F}(u)}{\pk{B_{1/2,(d-1)/2} > 1/c} \overline{F}(c u)}^{-2}\le  \frac{\overline{H}(u)}{\overline{Q}_d^2(u)} \le
\fracl{\overline{F}(u)}{\pk{B_{1/2,(d-1)/2} > 1/c} \overline{F}(c u)}^2. $$
The $O$-regular variation of $\overline{F}$ implies that $0 < b_1 \le \overline{H}(u)/\overline{Q}_d^d(u)\le b_2 < \IF$ for some constants $b_1,b_2$ and
for all $u$ large, hence \eqref{eq:vn1:2} yields that $\E{v_n},n\ge 1$ is a bounded sequence. If $d\ge 2$ the proof follows utilising \eqref{eq:vn1} and the bounds on the ratio $\overline{F}(u)/\overline{Q}_d(u)$,
and thus the result follows.
\QED

\prooftheo{theo:Hu} The proof can be established along the lines of Hueter (1999) utilising further \netheo{theo:m},
 the fact that $\xi_{Q_d}(n), Q_d^{-1}(1-1/n)$ are slowly regularly varying functions and \eqref{tn2}.  \QED

\section{Appendix}

{\bf A1}. First note that if $X^2 \sim Y^2 B_{1/2,1/2},$ with $Y >0$ almost surely independent of
$B_{1/2,1/2}$ and $X$ being symmetric about 0, then the distribution of $X$ is given by (see Carnal (1970))
$$ 2\pk{X> r} = \pk{\abs{X}> r}= \frac{2}{ \pi} \int_r^\IF \arccos(r/s)\, \abs{d \overline{G}(s)}, \quad r\ge 0,$$
with $G$ the distribution function of $Y$. Next, since $H$ is defined by (see Carnal (1970))
$$ H(r)= \frac{2}{\pi}\int_r^\IF \arccos(r/s)\,\abs{d \overline{F}^2(s)}, \quad r\ge 0$$
Eq. \eqref{A1} follows easily.

{\bf A2}. Let $R \sim H$ be a positive random variable independent of $B_{\alpha,\beta}, \alpha,\beta \in (0,\IF),$ and denote by
$Q_{\alpha, \beta}$ the distribution function of $ R B_{\alpha,\beta}$. Then $Q_{\alpha,\beta}$ possesses the density function
$q_{\alpha,\beta}$ given by  (see (22) in Hashorva et al. (2007))
\BQN \label{PKS:14:2}
q_{\alpha,\beta}(x)&=& \frac{\Gamma(\alpha+ \beta)}{\Gamma(\alpha)} x^{\alpha-1}  \int_x^\IF (s-x)^{\beta -1} s^{- \alpha -\beta+1} \, dH(s), \quad  \forall x\in (0,x_H).
\EQN
It is thus clear that  $Q_{\alpha, \beta}$ is a continuous distribution function.\\
Next assume that $N$ is a univariate distribution function with $N(0)=0$ and upper endpoint $x_N\in (0, \IF]$
such that $\mu_N:=\int_0^\IF x \, d N(x)\in (0,\IF)$. Define a new distribution function $N^*$ by
$$ N^*(s)= 1- \int_s^\IF \sqrt{x^2- s^2}\, d N(s)/\mu_N, \quad s\ge 0.$$
Then we have (see (18.5) in Reiss and Thomas (2007)) that $N^*$ possesses a density function $n^*$ given by
$$ n^*(s)=\frac{s}{\mu_N} \int_s^\IF(x^2- s^2)^{-1/2}\, d N(s), \quad s\ge 0$$
implying that $N^*$ is a continuous distribution function.

{\bf A3}. We give below two lemmas followed by two propositions, which are utilised in the proofs above. Note in passing that
the next lemma has been useful when dealing with the asymptotics of near extremes, see e.g.\ Pakes (2000). Furthermore, refinements under stronger asymptotic assumptions can be found in Li (2008).\\

\BL  \label{lem:2}
Let $G_1,G_2$ be two continuous distribution functions with upper endpoint $\omega\in (- \IF, \IF]$. Further let
$\lambda\in [0,1], l, \rho \in [0,\IF)$ be given constants and  $\kal{L}$ be a positive slowly varying function  at 0. Then the following statements are equivalent:\\

a) For any $z\inr$
\BQN\
\int_{\R} [ \lambda + (1- \lambda) \overline{G}_1(s)]^{n- z} \, dG_2(s) &\sim &
\frac{l}{((1- \lambda)n)^{\rho}} \kal{L}(1/n), \quad n\to \IF.
\EQN
b)
\BQN\label{lem2:2}
\frac{ \overline{G}_2(u)}{\overline{G}_1^\rho(u)\kal{L}(\overline{G}_1(u))}
&\sim &\frac{l }{\Gamma(\rho+1)}, \quad u \uparrow \omega.
\EQN
c)
\BQN
\int_{\R} \overline{G}_2(s) \exp(- (n-z) \overline{G}_1(s)) \, dG_1(s) &\sim &
 \frac{l}{n^{\rho-1}} \kal{L}(1/n), \quad n\to \IF.
\EQN
\EL
\prooflem{lem:2} The equivalence of the first two statements follows easily by Lemma 3.1 in Hashorva (2002). The
equivalence  of $b)$ and $c)$ can be established with similar arguments. \QED

\BL\label{lem:OR} Let $X \equaldis R Y$ with $R\sim F$ being independent of $Y\in [0,1]$ almost surely.
Suppose that $F(0)=0, x_F=\IF$ and $\pk{Y> s}\in (0,1)$ for any $s\in (0,1)$ and denote by $G$ the distribution function of $X$. Then $F$ is $O$-regularly varying iff $G$ is $O$-regularly varying.
\EL
\prooflem{lem:OR} By the assumptions we have for any $\alpha>1, x>0$
\BQNY
\overline{F}(x) \ge \overline{G}(x)\ge \int_{\alpha x} ^\IF \pk{Y > x/r}\, d F(r) \ge \pk{Y> 1/\alpha} \overline
{F}(\alpha x)>0
\EQNY
implying thus for any $c>1$
\BQNY
\frac{\overline{F}(cx)}{\pk{Y> 1/\alpha} \overline
{F}(\alpha x)}\ge  \frac{\overline{G}(cx)}{\overline{G}(x)} \ge  \pk{Y> 1/\alpha} \frac{\overline
{F}(\alpha cx)}{\overline{F}(x)}>0.
\EQNY
Choosing $\alpha\in (1,c)$ the assumption $F$ is $O$-regularly varying yields that
also $G$ is $O$-regularly varying. The proof of the converse follows with similar arguments, therefore it is omitted here.\\
 \QED


\BT\label{eq:theo:BM1} Let  $Y\sim H$ be independent of $B_{a,b}, a,b\in (0,\IF)$ with $H$ such that
$x_H \in (0, \IF],$ and $H(0)=0$. Let $X\equaldis Y [1-B_{a,b} ]^{1/\tau}, \tau\in (0,\IF)$ with distribution function $F$. Then we have:\\
i) $H\in GMDA(w) $ is equivalent with  $F\in GMDA(w) $, and
\BQN\label{res:BERM1:a}
\overline{F}(u) &\sim& \frac{\Gamma(a+b)}{\Gamma(b)} \fracl{\tau}{
uw(u)} ^{a}\overline{H}(u), \quad u \uparrow x_F.
\EQN
Furthermore, $F$ possesses a density function $f$ such that $f(u) \sim w(u) \overline{F}(u), u \uparrow x_F$.\

ii) The distribution function $H$ satisfies \eqref{eq:fre1} with some $\gamma \le 0$ iff
$F$ satisfies \eqref{eq:fre1} with the same $\gamma$, and moreover
\BQN
\label{res:BERM1:b}
\overline{F}(u) &\sim&
\frac{\Gamma(a+b)}{\Gamma(b)}\frac{\Gamma(b+\gamma/\tau)}{\Gamma(a+b+\gamma/\tau)}\overline{H}(u) , \quad u\to \IF.
\EQN
iii) The distribution $F$ with $x_F=1$ satisfies \eqref{eq:wei} with some $\gamma\ge 0$,
\COM{\BQN
\frac{ \overline{F}(1- su)}{\overline{F}(1- u)} \sim s^\gamma,  \quad \forall
s\in (0,\IF), \quad u\downarrow 0,
\EQN
with  $\gamma \le 0$ }
iff $H$ satisfies \eqref{eq:wei} with $\gamma^*:=\gamma+a,  \gamma\ge 0$, and moreover
\BQN\label{res:BERM1:c}
\overline{F}(u)
&\sim& \frac{\Gamma(a+b)}{\Gamma(b)} \frac{\Gamma(\gamma+1)}{\Gamma(\gamma+a+1)} (\tau (1- u))^a \overline{H}(u) ,  \quad
u \uparrow 1.
\EQN
\ET

\prooftheo{eq:theo:BM1} The proof follows from Theorem 16 in Hashorva et al.\ (2007) and the results of
Hashorva and Pakes (2010). \QED

\BT \label{theo:0}
Let $F, H,K, q_d, Q_d,d\ge 2$ be as in \netheo{theo:m}. Then  $F \in GMDA(w)$ iff one of the following relations hold:\\
a) For any $d\ge 2$ we have $Q_d \in GMDA(w)$. Furthermore
\BQN
q_d(u) &\sim &  w(u)  \overline{Q}_d(u), \quad u \uparrow x_F
\EQN
and
\BQN\label{theo:0:a}
\overline{Q}_d(u) 
&\sim & 2^{(d-3)/2} \frac{\Gamma(d/2)}{\sqrt{ \pi}}  ( u w(u))^{-(d-1)/2} \overline{F}(u), \quad u\uparrow x_F.
\EQN
b)  $H \in GMDA(2w)$ and moreover
\BQN\label{theo:0:b}
\overline{H}(u) &\sim & \frac{1}{\sqrt{\pi  u w(u)}} \overline{F}^2(u) \sim 2 \sqrt{\pi u w(u)}( \overline{Q}(u))^2, \quad u\uparrow x_F.
\EQN
c) We have $\mu_F:=\int_0^\IF y \, d F(y)$ is finite and $K^*(s)= 1- \int_s^\IF \sqrt{y^2- s^2} \, d F(y)/\mu_F, s\ge 0$ is a continuous distribution function with $K^*\in GMDA(w)$. Furthermore,  we have
\BQN\label{eq:carnal:miss}
\overline{K}(u) &\sim & \frac{1}{2} \overline{F}^2(u)  u/w(u) \sim   \pi u^2 \overline{Q}^2(u), \quad u\uparrow x_F.
\EQN
\ET
\prooftheo{theo:0} \COM{
Suppose that  $X \equaldis Y B,$ with $Y$ some positive random variable with upper endpoint $\omega \in \{1, \IF\}$ being independent of $B_{\alpha, \beta}$.  From Theorem 3.1 in Hashorva and Pakes (2010) $X$ has distribution function in the Gumbel MDA  with scaling function $w$ if  and only if $Y$ has distribution function in the Gumbel MDA  with the same scaling function $w$. Furthermore,
\BQNY
\pk{X> u}& \sim  \frac{\Gamma(\alpha+ \beta)}{\Gamma(\alpha)} 2^\beta (u w(u)) ^{-\alpha}\pk{Y> u}, \quad u\uparrow \omega.
\EQNY
}
 The proofs of statement $a)$ and $b)$ follow immediately by applying statement $i)$ of \netheo{eq:theo:BM1} in Appendix {\bf A3} (recall
that both \eqref{eq:camb} and \eqref{A1} hold).

Next, if $F \in GMDA(w)$,  then  $\mu_F=\int_0^\IF x\, d F(x)$ is finite, therefore $K^*(s)= 1- \int_s^\IF \sqrt{y^2- s^2} \, d F(y)/\mu_F, s\ge 0$ defines a distribution function with upper endpoint $x_F$.
Applying Lemma 7.6 in Hashorva (2007) we obtain
\BQNY
\overline{K}^*(u)&=& \frac{1}{\mu_F}\int_u^\IF \sqrt{x^2- u^2} \, d F(x) \sim \frac{1}{\mu_F} \Gamma(3/2) (2 u/w(u))^{1/2}\overline{F}(u), \quad
u\uparrow x_F.
\EQNY
Since locally uniformly in $\R$
\BQN\label{eq:loc}
\frac{w(u+ s/w(u))}{w(u)} &\to& 1, \quad u\uparrow x_F
\EQN
for any $s\inr,$ it follows that $K^*\in GMDA(w)$.\\
In order to finish the proof we need to show the converse. Assume therefore $K^*\in GMDA(w)$ and $\mu_F\in (0,\IF)$.  By the Abel integral equation
(see Heinrich (2007))
\BQNY
\overline{F}(x)&=& \frac{2\mu_F }{\pi} \int_x^\IF(y^2- x^2)^{-1/2} \, d K^*(x), \quad x\ge 0.
\EQNY
Applying again Lemma 7.6 in Hashorva (2007) we obtain
\BQNY
\overline{F}(u)&\sim& \frac{2\mu_F}{\pi} \Gamma(1/2) (2 u/w(u))^{-1/2}\overline{K}^*(u),
\quad  u\uparrow x_F,
\EQNY
hence the result follows. \QED

{\bf Acknowledgement}: I would like to thank the referee for a very kind and deep review and several suggestions which improved the presentation substantially. I am thankful to
Ilya Molchanov for many helpful discussions and comments, and to J\"urg H\"usler and  Michael Mayer  for providing some key references.

\bibliographystyle{plain}

\end{document}